\documentclass[fleqn]{elsart}

\journal{J. Comput. Appl. Math.}

\usepackage{amssymb,latexsym,amsmath,amsfonts}

\def\ds{\displaystyle}

\def\deg{\textrm{deg}}
\def\tr{\textrm{tr}}
\def\sign{\textrm{sign}}

\def \C{\mathbb{C}}
\def \R{\mathbb{R}}
\def \N{\mathbb{N}}

\def \reff#1{(\ref{#1})}
\def \v#1{\vec{#1}}

\numberwithin{equation}{section}
\newtheorem{theorem}{Theorem}[section]
\newtheorem{lemma}[theorem]{Lemma}
\newtheorem{Definition}[theorem]{Definition}
\newenvironment{definition}{\begin{Definition}\rm}{\end{Definition}}
\newtheorem{Remark}[theorem]{Remark}
\newenvironment{remark}{\begin{Remark}\rm}{\end{Remark}}
\newtheorem{Assumption}[theorem]{Assumption}

\newenvironment{proof}%
{\rm \trivlist \item[\hskip \labelsep{\bf Proof. }]}%
{\hspace*{\fill}$\Box$\endtrivlist}

\begin{document}
\begin{frontmatter}
\title{Gaussian quadrature for multiple orthogonal polynomials\thanksref{label1}}
\author{Jonathan Coussement\thanksref{label2}}
and
\ead{Jonathan.Coussement@wis.kuleuven.ac.be}
\author{Walter Van Assche}
\ead{Walter.VanAssche@wis.kuleuven.ac.be}
\thanks[label1]{This work was supported by INTAS Research Network 03-51-6631
and by FWO projects G.0184.02 and G.0455.04}
\thanks[label2]{J.\ Coussement is a research assistant of the Fund for Scientific Research --
Flanders (Belgium)}
\address{Katholieke Universiteit Leuven, Department of Mathematics,
Celestijnenlaan 200B, B-3001 Leuven, Belgium}
\begin{abstract}
    We study multiple orthogonal polynomials of type I and type II which
    have ortho\-gonality conditions with respect to $r$ measures.
    These polynomials are connected by their recurrence relation
    of order $r+1$.  First we show a relation with the eigenvalue
    problem of a banded lower Hessenberg matrix $L_n$, containing the
    recurrence coef\-ficients. As a consequence,
    we easily find that the multiple orthogonal polynomials of type I and
    type II satisfy a generalized Christoffel-Darboux identity.  Furthermore,
    we explain the notion of
    multiple Gaussian quadrature (for proper multi-indices), which is an
    extension of the theory of Gaussian quadrature for orthogonal polyno\-mials
    and was introduced by C.\ F.\ Borges.  In particular
    we show that the quadrature points and quadrature weights can be expressed
    in terms of the eigenvalue problem of $L_n$.
\end{abstract}

\end{frontmatter}

\section{Introduction}
\label{intro}

Multiple orthogonal polynomials arise naturally in the theory of
simultaneous rational approximation, in particular in the
Hermite-Pad\'e approximation of a system of $r\in \N$ (Markov and
Stieltjes) functions \cite{Bruin1,Bruin2,Mahler,Nikishin}.  They
are a generalization of orthogonal polynomials in the sense that
they satisfy orthogonality conditions with respect to $r \in
\mathbb{N}$ measures $\mu_1,\ldots ,\mu_r$ for which all the
moments exist.  In the literature one can find already a lot of
ex\-amples of such polynomials
\cite{aptekarev,aptbran,Coussement,Bernd,Els}.  Normally the
measures are taken to be positive. However, in this paper we study
{\em formal} multiple orthogonal polynomials which means that we
allow complex measures.

We will only consider multiple orthogonal polynomials with respect
to {\em proper multi-indices}. The proper multi-index
corresponding to $n\in \mathbb{N}_0=\N\cup \{0\}$ is
    \[\v \nu_n=(\underbrace{m+1,m+1,\ldots,m+1}_{s\ \textrm{times}}
    ,\underbrace{m,m,\ldots,m}_{r-s\ \textrm{times}}) \ \in \mathbb{N}_0^r,\]
where $n=mr+s$, $0 < s \le r$. There exist two types of multiple
orthogonal polynomials, type I and type II. Let
$\Gamma_1,\ldots,\Gamma_r$ be the supports of the $r$ measures. A
multiple orthogonal polynomial $P_{n}=P_{\v \nu_n}$ of type II
with respect to the proper multi-index $\v \nu_n$, is a polynomial
of degree at most $n$, which satisfies the orthogonality
conditions
   \begin{equation}
   \label{stelseltypeII}
   \int_{\Gamma_j} P_n(x)\:x^\ell \: d\mu_j (x) =0,
   \qquad \ell=0,\ldots ,\v \nu_n(j)-1,\quad j=1,\ldots ,r.
   \end{equation}
Here $\v \nu_n(j)$ is the $j$th component of $\v \nu_n$. Equation
(\ref{stelseltypeII}) leads to a system of $n$ homogeneous linear
equations for the $n+1$ unknown coefficients of $P_n$. A basic
requirement to have a good definition is that every possible
solution of the system (\ref{stelseltypeII}) has exactly degree
$n$.  This is equal to the assumption that \reff{stelseltypeII}
has a unique solution (up to a scalar multiplicative constant)
which has exactly degree $n$. In that case we call $\v \nu_n$ a
{\em normal multi-index} for $\mu_1,...,\mu_r$. Let
$m_\ell^{(j)}=\int_{\Gamma_j}x^\ell\: d\mu_j(x)$ be the $\ell$th
moment of the measure $\mu_j$ and set
    \begin{equation}
    \label{momentmatrix} D_{n} = \Bigl( D^{(1)}_{n,\v \nu_n(1)}\
    \cdots\  D^{(r)}_{n,\v \nu_n(r)} \Bigr)^T,
    \end{equation}
where
    \[
    D^{(j)}_{n,\ell} = \left(
    \begin{array}{cccc}
    m_0^{(j)} & m_1^{(j)} & \cdots & m_{\ell-1}^{(j)}\\
    m_1^{(j)} & m_2^{(j)} & \cdots & m_{\ell}^{(j)}\\
    \vdots & \vdots & & \vdots\\
    m_{n-1}^{(j)} & m_{n}^{(j)} & \cdots & m_{n+\ell-2}^{(j)}
    \end{array} \right)
    \]
is an $n\times \ell$ matrix of moments of the measure $\mu_j$.
Then $D_n$ is the matrix of the linear system
\reff{stelseltypeII}, without the last column.  It is known and
easily verified that $\v \nu_n$ is normal if and only if $D_n$ has
rank $n$ \cite{Coussement,Bernd,Nikishin}. In the case that all
the proper multi-indices $\v \nu_n$, $n\in \mathbb{N}$, are
normal, we call the system of measures a {\em weakly complete
system}.

A type I multiple orthogonal vector polynomial $\v A_n=\v A_{\v
\nu_n}= (A_{n,1},\ldots,A_{n,r})^T$, corres\-ponding to the proper
multi-index $\v \nu_n$, consists of $r$ polynomials $A_{n,j}$ of
degree at most $\v \nu_n(j)-1$, satisfying the orthogonality
conditions
    \begin{equation}
    \label{stelseltypeI}
    \int x^\ell \sum_{j=1}^{r}A_{n,j}(x)\ d\mu_j(x)=0,\qquad
    \ell=0,1,\ldots,n-2.
    \end{equation}
The usual requirement is that every $A_{n,j}$ has exactly degree
$\v \nu_n(j)-1$. However, as mentioned before, in this paper we
only consider proper multi-indices. So, if the system
\reff{stelseltypeI} has a unique solution up to a multiplicative
scalar constant and $A_{n,s}$ has exactly degree $m$, where
$n=mr+s$, $0 < s \le r$, we obtain a good definition. Note that
these conditions are satisfied if $\v \nu_{n-1}$ is normal, which
means that $D_{n-1}$ has rank $n-1$. This easily follows from the
fact that the matrix of the linear system \reff{stelseltypeI},
without the $s(m+1)$th column, is equal to $D_{n-1}^T$.

In Section \ref{sectionGauss} we approximate integrals of the kind
    \[\int_{\Gamma_j} f(x)\ d\mu_j(x),\qquad j=1,\ldots,r,\]
simultaneously, using weighted quadrature formulas which have the
same $n$ quadrature nodes.  C.\ F.\ Borges studied this already in
\cite{Borges}.  The point of interest is then to find the
appropriate quadrature nodes and weights maximizing the vector
order of the weighted quadrature formulas along the proper
multi-indices.  In \cite[Theorem 2]{Borges}, C.\ F.\ Borges showed
that this is obtained if we require each weighted quadrature
formula to be interpolating and the quadrature nodes to be the
zeros of the type II multiple orthogonal polynomial $P_n$, assumed
to be simple.  We recall this notion of {\em multiple Gaussian
quadrature} in Theorem \ref{maxorder}.

In the case of Gaussian quadrature ($r=1$) the nodes and weights
can be expressed in terms of the eigenvalue problem of a Jacobi
matrix, containing the recurrence coefficients, see e.g.\
\cite[Chapter 3, \S 2.3]{Gautschi} and \cite{Golub}. In Theorem
\ref{Gausseig} we extend this to multiple Gaussian quadrature.
Here we need some properties of multiple orthogonal polynomials.
In Section \ref{recursiemeervoudige} we recall that, for a weakly
complete system, the multiple orthogonal polynomials of type I and
type II each satisfy a recurrence relation of order $r+1$ which
are closely connected, see, e.g., \cite[\S 24]{Mahler} and
\cite{Brezinski,Kalyagin1,Sorokin,Iseghem}. As a consequence,
these polynomials of type I and II are then linked to the left and
right eigenvalue problem, respectively, of a banded lower
Hessenberg matrix $L_n$ containing the recurrence coefficients as
in \reff{L_n}.  This then also leads to a generalized
Christoffel-Darboux identity in Theorem \ref{Ch-D}, similar to
\cite[(21)]{Sorokin}. Using all this, we show that in the case of
multiple Gaussian quadrature the nodes are the eigenvalues of
$L_n$ and the weights can be expressed in terms of the
corresponding left and right eigenvectors.

\section{Recurrence relation of order $r+1$}
\label{recursiemeervoudige}

\subsection{Relation between type I and type II}

Suppose that the system of measures $\mu_1,\ldots,\mu_r$ forms a
weakly complete system.  The monic multiple orthogonal polynomials
$P_n$ of type II corresponding to proper multi-indices are then
uniquely determined.  Furthermore, as mentioned in the
introduction, also the multiple orthogonal polynomials of type I,
$\v A_n$, are unique up to a normalizing multiplicative constant.
In this paper we take the normalization
    \begin{equation}
    \label{normalisation_type_I}
    \int x^{n-1}\sum_{j=1}^{r}A_{n,j}(x)\
    d\mu_j(x)=1.
    \end{equation}
Note that the integral in (\ref{normalisation_type_I}) cannot be
zero since $\vec{\nu}_n$ is a normal multi-index.

It is known that multiple orthogonal polynomials satisfy a
recurrence relation of order $r+1$, see, e.g., \cite[\S
24]{Mahler} and \cite{Brezinski,Kalyagin1,Sorokin,Iseghem}.  For
the monic multiple orthogonal polynomials of type II this is
    \begin{equation}
    \label{recursie}
    xP_n(x) = P_{n+1}(x)+\sum_{j=0}^r
    a_{n,j}\ P_{n-j}(x),\qquad n\ge 0,
    \end{equation}
with initial conditions $P_0\equiv 1$ and $P_j\equiv 0$,
$j=-1,-2,\ldots,-r$.  The proof is similar as in the case of
orthogonal polynomials.  Furthermore, there exist some integral
representations for the recurrence coefficients in terms of the
multiple orthogonal polynomials of type I  and type II.  If we
integrate \reff{recursie} with respect to the measures
$\sum_{\ell=1}^r A_{n+1-k,\ell}(x)\ d\mu_\ell(x)$,
$k=0,1,\ldots,r$, then we find
        \begin{multline*}
        \int xP_n(x)\sum\limits_{\ell=1}^rA_{n+1-k,\ell}(x)\
        d\mu_\ell(x)
         =
        \int P_{n+1}(x)\sum\limits_{\ell=1}^{r}A_{n+1-k,\ell}(x)\
        d\mu_\ell(x)\\
        +\sum_{j=0}^r a_{n,j}\int P_{n-j}(x)\sum\limits_{\ell=1}^{r}
        A_{n+1-k,\ell}(x)\ d\mu_\ell(x).
        \end{multline*}
Now apply the orthogonality relations of the multiple orthogonal
polynomials of type I and II.  Because of the normalization
(\ref{normalisation_type_I}) and the fact that the $P_n$ are
monic, we then obtain
    \begin{equation}
    \label{integraltypeII}
    a_{n,j} =
    \int x P_n(x)\sum\limits_{\ell=1}^{r}A_{n+1-j,\ell}(x)\
    d\mu_\ell(x), \qquad j=0,1,\ldots,\min (r,n).
    \end{equation}

The multiple orthogonal polynomials of type I also satisfy a
recurrence relation of order $r+1$ which is related to the one of
type II.  Denote by  $\C^r[x]$ the space of $r$-dimensional vector
polynomials with complex coefficients.  The vector polynomials $\v
A_n$ then form a basis for $\C^r[x]$, so that there exist unique
$c_{n,l}$ for which $x\v A_n(x) = \sum_{k=1}^{n+r} c_{n,k}\ \v
A_k(x)$. By linearity we then have
        \[\int xP_i(x)\sum_{\ell=1}^{r}A_{n,\ell}(x)\ d\mu_\ell(x)
        = \sum_{k=1}^{n+r}
        c_{n,k}\int P_i(x)\sum_{\ell=1}^{r}A_{k,\ell}(x)\ d\mu_\ell(x),
        \qquad i=0,1,\ldots
        \]
The cases $i=0,1,\ldots,n-3$ give rise to
$c_{n,1}=\cdots=c_{n,n-2}=0$ by the orthogonality relations
\reff{stelseltypeI}.  Furthermore, use the orthogonality relations
\reff{stelseltypeII} and the normalization
(\ref{normalisation_type_I}) to get
    \[c_{n,n-1}=1,\qquad  c_{n,n+j}=\int xP_{n+j-1}(x)
    \sum_{\ell=1}^{r}A_{n,\ell}(x)\ d\mu_\ell(x), \qquad j=0,1,\ldots,r.\]
Comparing this with \reff{integraltypeII} we obtain
$c_{n,n+j}=a_{n+j-1,j}$, $j=0,1,\ldots,r$.  So if the monic
multiple orthogonal polynomials of type II satisfy the recurrence
relation \reff{recursie} then the multiple orthogonal polynomials
of type I with normalization (\ref{normalisation_type_I}) satisfy
the recurrence relation
        \begin{equation}
        \label{recursieI}
         x\v A_n(x)  =  \v A_{n-1}(x)+\sum\limits_{j=0}^r
        a_{n+j-1,j}\ \v A_{n+j}(x),\qquad n\ge 1,
        \end{equation}
with initial conditions $\v A_0\equiv \v 0$ and $\v A_1,\ldots,\v
A_r$. Denote by $D_n^{(j,i)}$ the matrix obtained from $D_n$ by
deleting the $j$th row and the $i$th column (and set $\det
D_1^{(1,1)}$ equal to $1$). For $1\le i,j\le r$ we then have
    \begin{equation}
    \label{initial}
    A_{i,j}=\left\{
    \begin{array}{lll}
    (-1)^{i+j}\ {\ds \frac{\det D_i^{(j,i)}}{\det D_i}}\:,  &  \qquad & j\le
    i,\\
    0\:, &  \qquad &  j>i,
    \end{array}
    \right.
    \end{equation}
which are functions of the (first $r$) moments of the measures
$\mu_1,\ldots ,\mu_r$.  Notice that $A_{j,j}\neq 0$,
$j=1,\ldots,r$, because the multi-indices $\v \nu_n$ are normal.
    \begin{remark}
    \label{nietnul}
    Let $n=mr+s$, $0<s\le r$.  Applying the orthogonality conditions
    of the polynomial $P_{n+r-1}$,  we obtain from \reff{integraltypeII}
        \[a_{n+r-1,r}  =
        \int x P_{n+r-1}(x)\sum\limits_{\ell=1}^{r}A_{n,\ell}(x)\ d\mu_\ell(x)\\
         =  \int x P_{n+r-1}(x)A_{n,s}(x)\ d\mu_s(x),\quad n\ge
        1.\]
    We assumed that $\v \nu_{n-1}$
    and $\v \nu_{n+r}$ are normal, so $A_{n,s}$ has exactly degree $m$
    and
        \begin{equation}
        a_{n+r-1,r}\not=0, \qquad n\ge 1.
        \end{equation}
    This means that we can recover the vector-polynomials $\v A_n$ by
    \reff{recursieI} if we know the initial conditions \reff{initial} and
    the recurrence coefficients.
    \end{remark}
    \begin{remark}
    \label{helesetvanmaten}
    We supposed that the measures $\mu_1,\ldots,\mu_r$ form a
    weakly complete system.  If we replace these measures by
        \begin{equation}
        \label{LCmeasures}
        \alpha_{1,1}\:\mu_1,\ \alpha_{2,1}\:\mu_1+\alpha_{2,2}\:\mu_2,
        \ \ldots,\ \sum_{j=1}^r\alpha_{r,j}\:\mu_j,
        \end{equation}
    where $\alpha_{i,j}\in \C$, $1\le j\le i\le r$, then it is clear that
    this system of measures forms a weakly complete system if
    and only if $\alpha_{j,j}\not=0$, $j=1,\ldots,r$.  Furthermore,
    all these systems of measures have the same multiple orthogonal
    polynomials of type II with proper multi-indices and so the same
     recurrence coefficients.  It is obvious that every
    set of measures of the form \reff{LCmeasures} then corresponds to
    another choice of initial conditions for the recurrence relation
    \reff{recursieI} where $A_{j,j}\neq 0$, $j=1,\ldots,r$.
    \end{remark}

\subsection{Eigenvalue problem of the banded Hessenberg matrix $L_n$}
\label{eigenvalsection}

In this section we introduce the banded lower Hessenberg matrix
    \begin{equation}
    \label{L_n}
    L_{n}=\left(
    \begin{array}{cccccccccc}
    a_{0,0} & 1 & 0 & \ldots & \ldots & \ldots & \ldots & \ldots & \ldots & 0\\
    a_{1,1} & a_{1,0}& 1 & 0 &  & & & & & 0\\
    a_{2,2} & a_{2,1}& a_{2,0} & 1 & 0 & & & & & 0\\
    \vdots &  & \ddots &  \ddots & \ddots & \ddots & & & & \vdots\\
    a_{r,r} & \ldots & \ldots & a_{r,1} & a_{r,0} & 1 & 0 & & & 0\\
    0 & a_{r+1,r} & \ldots & \ldots & a_{r+1,1} & a_{r+1,0} & 1 & 0 & & 0\\
    \vdots & \ddots & \ddots &  & & \ddots & \ddots & \ddots & \ddots & \vdots\\
    \vdots & & \ddots & \ddots &  & & \ddots & \ddots & \ddots & 0\\
    \vdots & & & \ddots & \ddots & & & \ddots & \ddots & 1\\
    0 & \ldots & \ldots & \ldots & 0 & a_{n-1,r} & \ldots & \ldots & a_{n-1,1} &
    a_{n-1,0}
    \end{array} \right),
    \end{equation}
which contains the recurrence coefficients of the multiple
orthogonal polynomials. If we expand the determinant
$\det(xI_n-L_n)$ along the last row, then we find that these
polynomials satisfy the recurrence relation \reff{recursie}. So
$P_n(x)= \det\left(xI_n-L_n\right)$ and the eigenvalues
$x_{1,n},\ldots,x_{n,n}$ of the matrix $L_n$ then coincide with
the zeros of the polynomial $P_n$.  We now have a closer look at
the relation between the type I and type II multiple orthogonal
polynomials and the eigenvectors of $L_n$.
The first $n$ relations of the recurrence \reff{recursie} can be
written as
    \begin{equation}
    \label{matrixrecursieII}
    L_n\ \left(\begin{array}{c}
           P_0(x)\\ \vdots\\ P_{n-2}(x)\\ P_{n-1}(x)
         \end{array}\right)
    +    \left(\begin{array}{c}
           0\\ \vdots\\ 0\\ P_{n}(x)
         \end{array}\right)
    =x\ \left(\begin{array}{c}
        P_0(x)\\ \vdots\\ P_{n-2}(x)\\ P_{n-1}(x)
        \end{array}\right).
    \end{equation}
Substituting $x=x_{\ell,n}$, $\ell=1,\ldots,n$, we then conclude
that
    $( P_0(x_{\ell,n})\ \cdots \
    P_{n-1}(x_{\ell,n}))^T$
is the (only possible) right eigenvector of $L_n$ corresponding to
the eigenvalue $x_{\ell,n}$, normalized so that the first
component is equal to $1$.  Note that the first component of a
right eigenvector of $L_n$ can not be $0$, otherwise each
component would be $0$.

For the type I polynomials we can find something similar. Define
for $ i=1,\ldots,r$ the vector polynomials
    \begin{align}
    \label{vecB1}
    \v B_n^{(i)}(x) & =  x\v A_{n-i+1}(x)-\v A_{n-i}(x)
    -\sum_{j=0}^{i-1} a_{n+j-i,j}\:\v A_{n-i+j+1}(x)\\
    \label{vecB2}
    & =  \sum_{j=i}^{r} a_{n-i+j,j}\:\v A_{n-i+j+1}(x),
    \end{align}
where we use definition \reff{vecB2}, with $a_{\ell_1,\ell_2}=1$,
$0\le \ell_1<\ell_2\le r$, in the cases $1\le n < i\le r$. The
first $n$ relations of the recurrence \reff{recursieI} then are
    \begin{multline}
    \label{matrixrecursieI}
    \Bigl(\v A_1(x)\ \cdots\ \v A_n(x)\Bigr)\ L_n
    +
    \Bigl(\underbrace{\v 0\ \cdots\ \v 0}_{n-r^\star} \ \v B_n^{(r^\star)}(x)\ \cdots\
    \v B_n^{(1)}(x)\Bigr)\\
    =x\ \Bigl(\v A_1(x)\ \cdots\ \v A_n(x)\Bigr),
    \end{multline}
where $r^\star=\min(r,n)$.
    \begin{definition}
    \label{qdef}
    Construct the polynomials $\v B_n^{(i)}$, $i=1,\ldots,r$, as in \reff{vecB1}
    and define for $n\in \N$ and $i=1,\ldots,\min(r,n)$ the polynomials
        \begin{equation}
        \label{Qveeltermen}
        Q^{(i)}_{k,n}(x)=\det \Bigl(\v B_n^{(1)}(x)\ \cdots\
        \v B_n^{(i-1)}(x)\ \v A_k(x)\ \v B_n^{(i+1)}(x)\ \cdots\
        \v B_n^{(r)}(x)\Bigr),
        \end{equation}
    $k\in \N$. We also define
        \begin{equation}
        B_n(x)=\det \Bigl(\v B_n^{(1)}(x)\ \cdots\
        \v B_n^{(r)}(x)\Bigr),\qquad n\in\N.
        \end{equation}
    \end{definition}
Note that each linear combination of the components of the vector
polynomials $\v A_k$ satisfies a similar relation as in
\reff{matrixrecursieI}.  In particular, for the polynomials
$Q^{(i)}_{k,n}$, $i=1,\ldots,\min(r,n)$, we get
    \begin{equation}
    \label{matrixrecursieQ}
    \Bigl(Q^{(i)}_{1,n}(x)\ \cdots\ Q^{(i)}_{n,n}(x)\Bigr)\ L_n
    +
    B_n(x)\ \v e_{n-i+1}^{\ T}
    =x\ \Bigl(Q^{(i)}_{1,n}(x)\ \cdots\ Q^{(i)}_{n,n}(x)\Bigr),
    \end{equation}
where $\v e_\ell$, $1\le \ell\le n$ is the $\ell$th unit vector in
$\R^n$. In Lemma \ref{verbandBenP} below we will prove that the
polynomial $B_n$ is equal to $P_n$ up to a multiplicative nonzero
constant. Equation \reff{matrixrecursieQ} then implies that the
vectors
    \begin{equation}
    \label{lefteig}
    \Bigl(Q^{(i)}_{1,n}(x_{\ell,n}) \ \cdots \
    Q^{(i)}_{n,n}(x_{\ell,n})\Bigr)^T,\qquad i=1,\ldots,\min(r,n),
    \end{equation}
are left eigenvectors of the matrix $L_n$ corresponding to the
eigenvalue $x_{\ell,n}$, if they are different from $\vec{0}$.
    \begin{lemma}
    \label{verbandBenP}
    Let $\v A_n$ the vector polynomials defined by the recurrence relation
    (\ref{recursieI}) and its initial conditions (with $A_{j,j}\neq 0$,
    $j=1,\ldots,r$).  Here we assume that the recurrence coefficients
    $a_{\ell,r}$, $\ell\ge r$, are different from 0.
    For the polynomials $B_n$ we then have that $B_n(x)=\gamma_n\:P_n(x)$,
    where the $P_n$ are the monic polynomials
    defined by the recurrence relation
    \reff{recursie} and
        \begin{equation}
        \label{gamma}
        \gamma_n=
        \frac{\prod_{j=1}^{r}A_{j,j}}{\prod_{\ell=r}^{n-1}a_{\ell,r}}
        \ (-1)^{\lfloor \frac{s}{2} \rfloor+\lfloor
        \frac{r-s}{2}\rfloor},
        \qquad n=mr+s,\quad 0<s\le r.
        \end{equation}
    \end{lemma}
    \begin{proof}
    First of all we prove that $B_n$ is a polynomial
    of exactly degree $n$ with leading coefficient $\gamma_n$.
    Using \reff{vecB2} we find that
        \begin{align}
        \nonumber
        B_n(x) & =  \det \Bigl(\v B_n^{(1)}(x)\ \cdots\
                        \v B_n^{(r)}(x)\Bigr)\\
        \label{dedet}
                & =  \prod_{\ell=n}^{n+r-1} a_{\ell,r}\ \det \Bigl(\v A_{n+r}(x)\ \cdots\
                        \v A_{n+1}(x)\Bigr)\\
        \label{expdet}
                & =  \prod_{\ell=n}^{n+r-1} a_{\ell,r}\  \sum_{\pi
                \in S_r} \sign (\pi)\
                A_{n+r,\pi(1)}(x)\:\ldots\:A_{n+1,\pi(r)}(x),
        \end{align}
    where we denote by $S_r$ the set of permutations of $r$
    elements.  Note that from the recurrence relation (\ref{recursieI})
    and its initial conditions, with $A_{j,j}\neq 0$,
    $j=1,\ldots,r$ and $a_{\ell,r}\not= 0$, $\ell\ge r$, we obtain
    that
        \begin{equation}
        A_{mr+s,s}(x)=\frac{A_{s,s}}{\prod_{i=1}^m
        a_{ir+s-1,r}}\:x^m
        +\mathcal{O}\left(x^{m-1}\right),
        \qquad 1\le s\le r,\quad m\in\mathbb{N}\cup \{0\}.
        \end{equation}
    Combining this and \reff{expdet} we see that, with $n=mr+s$, $0<s\le
    r$,
        \[\deg (B_n(x))
        =\deg (\underbrace{A_{n+r,s}(x)\:\ldots\:A_{n+r-s+1,1}(x)}_{s}\:
        \underbrace{A_{n+r-s,r}(x)\:\ldots\:A_{n+1,s+1}(x)}_{r-s})=n,\]
    where we note that all the other terms in \reff{expdet} have lower
    degree.  The sign of the permutation corresponding to this term is $(-1)^{\lfloor \frac{s}{2} \rfloor+\lfloor
    \frac{r-s}{2}\rfloor}$ so that $\gamma_n$ is the leading coefficient of $B_n$.

    To complete the proof we show that the polynomials $\frac{B_n}{\gamma_n}$ satisfy
    the recurrence relation \reff{recursie}. Note that
        \begin{equation}
        \label{verhoudinggamma}
        \frac{\gamma_n}{\gamma_{n-j}}=\frac{(-1)^{(r-1)j}}{\prod_{\ell=n-j}^{n-1}\:a_{\ell,r}},\qquad
        1\le j \le n,
        \end{equation}
    where we assume $a_{0,r}=\cdots=a_{r-1,r}=1$.
    By \reff{verhoudinggamma} and \reff{dedet}, we then obtain that
        \begin{align}
        \nonumber & \gamma_n\left(\frac{B_{n+1}(x)}{\gamma_{n+1}}+(a_{n,0}-x)\frac{B_n(x)}{\gamma_n}
        +\sum_{j=1}^r a_{n,j}\:\frac{B_{n-j}(x)}{\gamma_{n-j}}\right)\\
        \nonumber & =
        (-1)^{r-1}\:a_{n,r}\:B_{n+1}(x)+(a_{n,0}-x)B_n(x)+\sum_{j=1}^r
        (-1)^{(r-1)j}\:a_{n,j}\:\frac{B_{n-j}(x)}{\prod_{\ell=n-j}^{n-1}\:a_{\ell,r}}\\
        \nonumber & =
        (-1)^{r-1}\:\prod_{\ell=n}^{n+r} a_{\ell,r}\:\det \Bigl(\v
        A_{n+r+1}(x)\:\cdots\:\v A_{n+2}(x)\Bigr)\\
        \nonumber & \qquad  +(a_{n,0}-x)\prod_{\ell=n}^{n+r-1} a_{\ell,r}\:\det \Bigl(\v
        A_{n+r}(x)\:\cdots\:\v A_{n+1}(x)\Bigr)\\
        \label{uitgeschreven} & \qquad  +\sum_{j=1}^r
        (-1)^{(r-1)j}\:a_{n,j}\:\prod_{\ell=n}^{n-j+r-1} a_{\ell,r}\:\det \Bigl(\v
        A_{n-j+r}(x)\:\cdots\:\v A_{n-j+1}(x)\Bigr).
        \end{align}
    Form the recurrence relation \reff{recursieI} we first of all
    get that
        \begin{align*}
        & (-1)^{r-1}\:a_{n+r,r}\:\det \Bigl(\v A_{n+r+1}(x)\:\cdots\:\v
        A_{n+2}(x)\Bigr)\\
        & =
        (-1)^{r-1}\:\det \Bigl(-\v A_n(x)-(a_{n,0}-x)\v A_{n+1}(x)\ \v A_{n+r}(x)\:\cdots\:\v A_{n+2}(x)\Bigr)\\
        & =
        -\det \Bigl(\v A_{n+r}(x)\:\cdots\:\v A_{n+2}(x)\:\v A_n(x)\Bigr)-
        (a_{n,0}-x)\det \Bigl(\v A_{n+r}(x)\:\cdots\:\v
        A_{n+1}(x)\Bigr),
        \end{align*}
    so that the first two terms of \reff{uitgeschreven} reduce to
    $-\prod_{\ell=n}^{n+r-1}a_{\ell,r}\:\det (\v A_{n+r}(x)\:\cdots\:\v A_{n+2}(x)\:\v
    A_n(x))$.  We apply this argument several times, knowing that, for $j=1,\ldots,r-1$,
        \begin{align*}
        & a_{n+r,r}\:\det \Bigl(\v A_{n+r+1-j}(x)\:\cdots\:\v
        A_{n+2}(x)\:\v A_n(x)\:\cdots\:\v A_{n-j+1}(x)\Bigr)\\
        & =
        (-1)^{r}\:\det \Bigl(\v A_{n+r-j}(x)\:\cdots\:\v
        A_{n+2}(x)\:\v A_n(x)\:\cdots\:\v A_{n-j}(x)\Bigr)\\
        & \qquad
        +(-1)^{r-j}\:a_{n,j}\:\det \Bigl(\v A_{n+r-j}(x)\:\cdots\:\v
        A_{n-j+1}(x)\Bigr).
        \end{align*}
    Finally we obtain that the expression in
    \reff{uitgeschreven} is equal to 0.  This proves the lemma.
    \end{proof}

\subsection{Generalized Christoffel-Darboux identity}

For orthogonal polynomials it is known that they satisfy the
Christoffel-Darboux formula, see, e.g., \cite{Chihara}. In
\cite{Sorokin} the authors proved a generalized
Christoffel-Darboux identity for matrix orthogonality of vector
polynomials.  This includes a generalized Christoffel-Darboux
identity for the multiple orthogonal polynomials of type I and
type II in the sense that one of the vector polynomials then has
just one component. Here we show that this identity can be found
as a natural consequence of Section \ref{eigenvalsection}.
    \begin{theorem}
    \label{Ch-D}
    Suppose that the measures $\mu_1,\ldots,\mu_r$ form a
    weakly complete system. For the corresponding multiple orthogonal
    polynomials we have
        \begin{equation}
        \label{Ch-Darboux}
        (x-y)\sum_{k=1}^nP_{k-1}(x)Q^{(i)}_{k,n}(y)=\gamma_n
        \Bigl(P_n(x)P_{n-i}(y)
        -P_{n}(y)P_{n-i}(x)\Bigr)
        \end{equation}
    and
        \begin{equation}
        \label{Ch-Darbouxy=x}
        \sum_{k=1}^nP_{k-1}(x)Q^{(i)}_{k,n}(x)=\gamma_n
        \Bigl(P'_n(x)P_{n-i}(x)-P_n(x)P'_{n-i}(x)\Bigr),
        \end{equation}
    $i=1,\ldots,\min(r,n)$, with $Q^{(i)}_{k,n}$ as in Definition \ref{qdef} and
    $\gamma_n$ as in \reff{gamma}.
    \end{theorem}
    \begin{remark}
    The Christoffel-Darboux formulas in
    Theorem \ref{Ch-D} are related to the ge\-neralized
    Christoffel-Darboux identity \cite[(21)]{Sorokin},
    with one vector polynomial having just 1 component. In
    particular, for each $i$, expression \reff{Ch-Darboux} is a
    linear combination of the vector components in
    \cite[(21)]{Sorokin}.
    \end{remark}
    \begin{proof}
    Define $\v P(x)=( P_0(x)\ \cdots \ P_{n-1}(x))^T$ and
    $\v Q^{(i)}(x)=( Q^{(i)}_{1,n}(x)\ \cdots \
        Q^{(i)}_{n,n}(x))^T$, $i=1,\ldots,\min(r,n)$.
    The equations (\ref{matrixrecursieII}) and (\ref{matrixrecursieQ})
    then reduce to
        \begin{equation}
        \label{matrixvectorrecursieII}
        L_n\v P(x)+P_n(x)\vec{e}_n=x\v P(x)
        \end{equation}
    and
        \begin{equation}
        \label{matrixvectorrecursieI}
        \v Q^{(i)}(x)^TL_n+B_n(x)\vec{e}_{n-i+1}^{\ T}=x\v Q^{(i)}(x)^T.
        \end{equation}
    Now multiply (\ref{matrixvectorrecursieII}) on the right by the vector
    $\v Q^{(i)}(y)^T$ and (\ref{matrixvectorrecursieI}), evaluated at $y$, on
    the left by $\v P(x)$, so that
        \begin{align*}
        L_n\v P(x)\v Q^{(i)}(y)^T + P_n(x)\vec{e}_n\v Q^{(i)}(y)^T & =  x \v P(x)\v
        Q^{(i)}(y)^T,\\[1ex]
        \v P(x)\v Q^{(i)}(y)^TL_n + B_n(y)\v P(x)\vec{e}_{n-i+1}^{\ T} & =  y \v P(x)\v Q^{(i)}(y)^T.
        \end{align*}
    Next, take the trace of these two equations and subtract.  Since for two arbitrary squared matrices $C$ and $D$
    the property $\tr (CD)=\tr (DC)$ holds, we then obtain
        \begin{equation}
        \label{mettrace}
        (x-y)\ \tr \left(\v P(x)\v Q^{(i)}(y)^T\right)=P_n(x)Q^{(i)}_{n,n}(y)
        -B_n(y)P_{n-i}(x).
        \end{equation}
    First of all we note that $B_n(y)=\gamma_n\:P_n(y)$ by Lemma \ref{verbandBenP}.  Secondly, using \reff{vecB1}
    and \reff{vecB2}, we get
        \begin{align*}
        & Q_{n,n}^{(i)}(y)\\
        & =
        \det \Bigl(\v B_n^{(1)}(y)\ \cdots\
        \v B_n^{(i-1)}(y)\ \v A_n(y)\ \v B_n^{(i+1)}(y)\ \cdots\
        \v B_n^{(r)}(y)\Bigr)\\
        & =
        \det \Bigl(\underbrace{-\v A_{n-1}(y)\ \cdots\
        -\v A_{n-i+1}(y)}_{i-1}\ \v A_n(y)\ \underbrace{a_{n+r-i-1,r}\v A_{n+r-i}(y)\ \cdots\
        a_{n,r}\v A_{n+1}(y)}_{r-i}\Bigr)\\
        & =  \prod_{\ell=n}^{n+r-i-1} a_{\ell,r}\:(-1)^{(r-1)i}\:
        \det \Bigl(\v A_{n-i+r}(y)\ \cdots\ \v
        A_{n-i+1}(y)\Bigr)\\
        & =
        (-1)^{(r-1)i}\:\left(\prod_{\ell=n-i}^{n-1}a_{\ell,r}\right)^{-1}\:B_{n-i}(y).
        \end{align*}
    Combining this with Lemma \ref{verbandBenP} and \reff{verhoudinggamma} we can conclude that
    $Q_{n,n}^{(i)}(y)=\gamma_n\:P_{n-i}(y)$.  So, from
    \reff{mettrace} we finally obtain expression \reff{Ch-Darboux}.
    If we let $y\to x$ we also find (\ref{Ch-Darbouxy=x}).
    \end{proof}

\section{Multiple Gaussian quadrature}
\label{sectionGauss}

Suppose $\mu_1,\ldots,\mu_r$ are measures for which all the
moments exist. In this section we want to approximate a set of
integrals of the kind
    \[\int_{\Gamma_j} f_j(x)\ d\mu_j(x),\qquad j=1,\ldots,r.\]
An obvious possibility is to approximate each of these integrals
separately using Gaussian quadrature.   However, when $f_1= \cdots
= f_r=f$ we can consider to take the same set of quadrature nodes
in each weighted quadrature formula.  Our approximation then looks
like
    \begin{equation}
    \label{weightedformula}
    \int_{\Gamma_j} f(x)\:d\mu_j(x)  =  \sum_{\ell=1}^n w^{(j)}_{\ell,n}\:f(x_{\ell,n})+E^{(j)}_{n}(f),
    \qquad j=1,\ldots,r,
    \end{equation}
for different quadrature nodes $x_{1,n},\ldots,x_{n,n}$, which is
due to C.\ F.\ Borges \cite{Borges}. So, we will not necessarily
have the maximal order for each of the weighted quadrature
formulas. However, the advantage of this choice is that we only
need $n$ evaluations of the function $f$ instead of $rn$. In
accordance with the case $r=1$ we say that the set of weighted
quadrature formulas (\ref{weightedformula}) has {\em vector order}
$\v d=(d_1,\ldots,d_r)\in \N_0^r$ if $E^{(j)}_{n}(f)=0$ for each
polynomial $f$ of degree less than or equal to $d_j$,
$j=1,\ldots,r$.

Our goal is to maximize the vector order of \reff{weightedformula}
along the proper multi-indices.  A necessary requirement is then
that each of the quadrature formulas is interpolating, which means
that we have vector order at least $(n-1,\ldots,n-1)$. This
corresponds to the conditions
    \begin{equation}
    \label{eiginterpol}
    w^{(j)}_{\ell,n}  =  \int_{\Gamma_j} l_{\ell,n}(x)\ d\mu_j(x),
    \qquad \ell=1,\ldots,n,\qquad j=1,\ldots,r,
    \end{equation}
where we denote by
    \[l_{\ell,n}(x)=\prod_{i=1,i\not= \ell}^n
    \frac{x-x_{i,n}}{x_{\ell,n}-x_{i,n}},\qquad \ell=1,\ldots,n,\]
the fundamental polynomials of Lagrange interpolation.  The only
freedom we still have then consists of the choice of the set of
quadrature nodes.  In \cite{Borges} C.\ F.\ Borges already proved
that the maximal vector order (along the proper multi-indices) is
obtained if we choose the quadrature nodes to be the zeros of the
type II multiple orthogonal polynomial $P_n$, corresponding to the
measures $\mu_1,\ldots,\mu_r$.  We recall this theorem and give a
proof for completeness. Of course we need that the proper
multi-index $\v\nu_n$ is normal so that $P_n$ is uniquely defined
and has exact degree. Furthermore, we need that the zeros of $P_n$
are simple. Note that, in the case of positive measures, these
conditions are satisfied if the measures form, for example, an AT
system or Angelesco system \cite{Nikishin,Els}.
    \begin{theorem}[Borges]
    \label{maxorder}
    Suppose that the measures $\mu_1,\ldots,\mu_r$ form a
    weakly complete system and that $P_n$, the type
    II multiple orthogonal polynomial of degree $n$, has simple
    zeros. We then speak of multiple Gaussian quadrature if
    the weighted quadrature formulas in \reff{weightedformula} are interpolating and the
    quadrature nodes are the zeros of $P_n$.
    In this case we have vector order $(n-1)\v e+\v \nu_n$, where $\v e=(1,\ldots,1)\in \N^r$.
    Furthermore, these quadrature formulas do not have vector order $(n-1)\v e+\v \nu_{n+1}$.
    \end{theorem}
    \begin{proof}
    We prove that, with the conditions of the theorem, the $j$th weighted quadrature formula in
    \reff{weightedformula} has order $n-1+\v \nu_n(j)$, $j=1,\ldots,r$.
    Let $h_j$ be a polynomial of degree at least $n$ and at most $n-1+\v
    \nu_n(j)$. (If $\v \nu_n(j)=0$, there is nothing to prove.)
    Denote by $T_{n-1}^{(j)}$ the interpolating polynomial of $h_j$
    at the zeros $x_{1,n},\ldots,x_{n,n}$ of the type II multiple orthogonal
    polynomial $P_n$.  The polynomial $h_j-T_{n-1}^{(j)}$ can then
    be written as
        \[h_j(x)-T_{n-1}^{(j)}(x)=P_n(x)R^{(j)}(x),\]
    where $R^{(j)}$ is a polynomial of degree at most
    $\v \nu_n(j) -1$. Since the $j$th weighted quadrature formula is interpolating,
    we have
        \begin{align*}
        \int h_j(x)\:d\mu_j(x) - \sum_{\ell=1}^n h_j(x_{\ell,n})\:w_{\ell,n}^{(j)}
        & = \int h_j(x)\:d\mu_j(x) -
        \sum_{\ell=1}^n T_{n-1}^{(j)}(x_{\ell,n})\:w_{\ell,n}^{(j)}\\
        & = \int \left(h_j(x)-T_{n-1}^{(j)}(x)\right)\:d\mu_j(x)\\
        & = \int P_n(x)R^{(j)}(x)\:d\mu_j(x).
        \end{align*}
    By the orthogonality conditions of the type II polynomial
    $P_n$ we then see that this is equal to 0.

    The $x_{\ell,n}$ are the zeros of the polynomial $P_n$, so
        \[\sum_{\ell=1}^n P_n(x_{\ell,n}) \sum_{j=1}^r
        A_{n+1,j}(x_{\ell,n})\:w_{\ell,n}^{(j)}=0.\]
    Since $\v \nu_{n+1}$ is a normal index, this is different from
    $\int P_n(x)\sum_{j=1}^r A_{n+1,j}(x)\:d\mu_j(x)$, which means
    that we do not have vector order $(n-1)\v e+\v \nu_{n+1}$.
    \end{proof}

We now study the theory of multiple Gaussian quadrature from the
practical point of view.  In particular we give a link with the
eigenvalue problem of the banded lower Hessenberg matrix $L_n$,
see \reff{L_n}. First of all, from Section \ref{eigenvalsection}
we know that the zeros of the polynomial $P_n$ are the eigenvalues
of $L_n$. In the theorem below we show that, for multiple Gaussian
quadrature, the quadrature weights can be expressed in terms of
the corresponding left and right eigenvectors.  This extends the
expressions found in the case of Gaussian quadrature, see e.g.\
\cite[Chapter 3, \S 2.3]{Gautschi} and \cite{Golub}.
    \begin{theorem}
    \label{Gausseig}
    Assume that the measures $\mu_1,\ldots,\mu_r$ form a
    weakly complete system and that the type
    II multiple orthogonal polynomial $P_n$ has simple
    zeros $x_{1,n},\ldots,x_{n,n}$.  Let $\vec{v}_{\ell,n}$ be the right
    eigenvector of the matrix $L_n$ corresponding to the eigenvalue $x_{\ell,n}$,
    with first component equal to 1.  Similarly, denote by $\v u_{\ell,n}$ the
    corresponding left eigenvector with the $k_\ell$th component, the first non-zero
    component, equal to 1.  Here $1\le k_\ell \le \min(r,n)$.  In the case of multiple Gaussian quadrature
    we then have, for $\ell=1,\ldots,n$,
        \begin{equation}
        \label{uitdrgauss}
        w^{(j)}_{\ell,n}=\frac{1}{\v u_{\ell,n}^{\:T}\v v_{\ell,n}}\:
        \left(\sum_{k=1}^{\min (j,n)}C_{j,k}\:\v u_{\ell,n}(k)\right),\qquad j=1,\ldots, r,
        \end{equation}
    where
        \begin{align}
        C_{j,k}
        & =
        \int_{\Gamma_j} P_{k-1}(x)\:d\mu_j(x)\\
        \nonumber & =  \sum_{i=1}^k
        (-1)^{k+i}\:m_{i-1}^{(j)}\:\frac{\det D_k^{(k,i)}}{\det
        D_{k-1}}, \qquad 1\le k\le j\le r.
        \end{align}
    Here the constants $C_{j,j}$, $j=1,\ldots,r$ are different from
    $0$ and $w^{(j)}_{\ell,n}=0$, $j=1,\ldots, k_\ell-1$.
    \end{theorem}
    \begin{proof}
    In this proof we fix $\ell \in \{1,\ldots,n\}$.  From Section
    \ref{eigenvalsection} we know that
        \[\v v_{\ell,n}=\Bigl( P_0(x_{\ell,n})\ \cdots \
        P_{n-1}(x_{\ell,n})\Bigr)^T.\]
    By Remark \ref{nietnul}
    we have $a_{m+r-1,r}\not=0$, $m\ge 1$.  Then it is clear that at least one of the first
    $\min(r,n)$ components of the left eigenvector $\v u_{\ell,n}$ and at least one of the
    last $\min(r,n)$ components of $\v v_{\ell,n}$  is different from 0.  So,
    there exists an $i_\ell \in \{1,\ldots,\min(r,n)\}$ for which
    $P_{n-i_\ell}(x_{l,n})\not= 0$.  Taking $x=x_{\ell,n}$ in \reff{Ch-Darbouxy=x} with $i=i_\ell$, we then find
        \begin{equation}
        \label{Ch-Darbouxx=nulpunt}
        \sum_{k=1}^nP_{k-1}(x_{\ell,n})Q^{(i_\ell)}_{k,n}(x_{\ell,n})=\gamma_n
        P'_n(x_{\ell,n})P_{n-i_\ell}(x_{\ell,n})\not=0,
        \end{equation}
    because $x_{\ell,n}$ is a simple zero of $P_n$.  Consequently, the vector
    \reff{lefteig} with $i=i_\ell$ is different from $\v 0$
    and
        \begin{equation}
        \label{lefteigexact}
        \v u_{\ell,n}=\Bigl(\underbrace{0\:\cdots\:0}_{k_\ell-1}\ 1\
        \frac{Q^{(i_\ell)}_{k_\ell+1,n}(x_{\ell,n})}
        {Q^{(i_\ell)}_{k_\ell,n}(x_{\ell,n})} \ \cdots \
        \frac{Q^{(i_\ell)}_{n,n}(x_{\ell,n})}
        {Q^{(i_\ell)}_{k_\ell,n}(x_{\ell,n})}\Bigr)^T.
        \end{equation}
    Applying \reff{Ch-Darbouxx=nulpunt} we then find that
        \begin{equation}
        \label{inproduct}
        \v u_{\ell,n}^{\:T}\v v_{\ell,n}=\frac{\gamma_n
        P'_n(x_{\ell,n})P_{n-i_\ell}(x_{\ell,n})}
        {Q^{(i_\ell)}_{k_\ell,n}(x_{\ell,n})}.
        \end{equation}
    Next, since the weighted quadrature formulas are assumed to be interpolating, we obtain from
    the Christoffel-Darboux formula \reff{Ch-Darboux} with
    $i=i_\ell$ and $y=x_{\ell,n}$ that
        \begin{align*}
        w^{(j)}_{\ell,n}
        & =
        \int_{\Gamma_j} l_{\ell,n}(x)\ d\mu_j(x)\\
        & =
        \frac{1}{P_n'(x_{\ell,n})}\int_{\Gamma_j}
        \frac{P_n(x)}{x-x_{\ell,n}}\:d\mu_j(x)\\
        & =
        \frac{1}{\gamma_n P_n'(x_{\ell,n})P_{n-i_\ell}(x_{\ell,n})}
        \sum_{k=1}^nQ^{(i_\ell)}_{k,n}(x_{\ell,n})\int_{\Gamma_j}
        P_{k-1}(x)\:d\mu_j(x),
        \qquad j=1,\ldots,r.
        \end{align*}
    Combining this with \reff{inproduct} and using the orthogonality
    conditions of the type II polynomial $P_n$ we then finally get
        \[w^{(j)}_{\ell,n}=\frac{1}{\v u_{\ell,n}^{\:T}\v v_{\ell,n}}
        \sum_{k=1}^{\min(j,n)}\frac{Q^{(i_\ell)}_{k,n}(x_{\ell,n})}
        {Q^{(i_\ell)}_{k_\ell,n}(x_{\ell,n})}\int_{\Gamma_j}
        P_{k-1}(x)\:d\mu_j(x),
        \qquad j=1,\ldots,r,\]
    which proves \reff{uitdrgauss}.  Finally we note that the proper multi-indices are normal,
    so we have that $C_{j,j}\not=0$, $j=1,\ldots,r$.
    \end{proof}
    \begin{remark}
    \label{one to one}
    There exists a one to one mapping between the sets of constants
    $C_{j,k}\in \C$, $1\le k\le j\le r$, (with $C_{j,j}\not=0$, $1\le j\le r$) and
    $A_{i,j}\in\C$, $1\le j\le i\le r$, (with $A_{j,j}\not=0$, $1\le j\le r$).
    By the orthogonality conditions of the multiple
    orthogonal polynomials of type I and type II we obtain that
        \[
        \sum_{j=k}^i A_{i,j}\:C_{j,k}=
        \int P_{k-1}(x)\:\sum_{j=1}^i A_{i,j}\:d\mu_j(x)=
        \left\{
        \begin{array}{lll}
        0\:, & \ & k=1,\ldots,i-1,\\
        1\:, & \ & k=i,
        \end{array}
        \right.
        \qquad i=1,\ldots,r.
        \]
    In particular, each of these sets can be found from the
    other one by solving the linear systems
        \begin{equation}
        \label{CnaarA}
        \left(
        \begin{array}{cccc}
            C_{1,1} & C_{2,1} & \cdots & C_{i,1}\\
            0 & C_{2,2} & \cdots & C_{i,2}\\
            \vdots & \ddots & \ddots & \vdots\\
            0 & \cdots & 0 & C_{i,i}
        \end{array}
        \right)
        \left(
        \begin{array}{c}
            A_{i,1}\\
            A_{i,2}\\
            \vdots\\
            A_{i,i}
        \end{array}
        \right)
        =\left(
        \begin{array}{c}
            0\\
            \vdots\\
            0\\
            1
        \end{array}
        \right),
        \qquad i=1,\ldots,r,
        \end{equation}
    and
        \begin{equation}
        \label{AnaarC}
        \left(
        \begin{array}{cccc}
            A_{i,i} & 0 & \cdots & 0\\
            \vdots & \ddots & \ddots & \vdots\\
            \vdots &  & \ddots & 0\\
            A_{r,i} & \cdots & \cdots & A_{r,r}
        \end{array}
        \right)
        \left(
        \begin{array}{c}
            C_{i,i}\\
            C_{i+1,i}\\
            \vdots\\
            C_{r,i}
        \end{array}
        \right)
        =\left(
        \begin{array}{c}
            1\\
            0\\
            \vdots\\
            0
            \end{array}
        \right),
        \qquad i=1,\ldots,r,
        \end{equation}
    respectively.
    \end{remark}

\section{Conclusion}
\label{conclusionsection}

In Theorem \ref{Gausseig} we proved that in the case of multiple
Gaussian quadrature the quadrature nodes and weights in
\reff{weightedformula} can be found explicitly by solving the left
and right eigenvalue problem of a banded lower Hessenberg matrix
$L_n$.  Here $L_n$ contains the recurrence coefficients of the
multiple orthogonal polynomials corresponding to the set of
measures $\mu_1,\ldots,\mu_r$ as in \reff{L_n}.  Further, we only
need a set of constants $C_{j,k}\in \C$, $1\le k\le j\le r$, (with
$C_{j,j}\not=0$, $1\le j\le r$) which can be found from the
initial values of the type I multiple orthogonal polynomials by
solving the linear systems \reff{CnaarA}.

By Theorem \ref{maxorder} the weighted quadrature formulas
preserve the orthogonality conditions (and normalization) for the
finite set of polynomials $\v A_1,\ldots,\v A_n$ and
$P_0,P_1,\ldots,P_{n-1}$.  So, these are also the type I and type
II multiple orthogonal polynomials corresponding to the set of
discrete measures
\[\mu_{j,n}=\sum_{\ell=1}^nw_{\ell,n}^{(j)}\delta_{x_{\ell,n}},\qquad j=1,\ldots,r,\]
with finite support.  Then note that Remark \ref{one to one} is a
nice illustration of Remark \ref{helesetvanmaten}.  Since we
explained how to find these discrete measures starting from the
recurrence coefficients and the set of initial values
$A_{i,j}\in\C$, $1\le j\le i\le r$, (with $A_{j,j}\not=0$, $1\le
j\le r$), this gives rise to a constructive proof of the spectral
theorem for a finite set of multiple orthogonal polynomials. In
the case of an infinite set some results for the spectral theorem
were already obtained in \cite{Kalyagin2,Sorokin}.  However,
finding necessary and sufficient conditions on the recurrence
coefficients to have positive orthogonality measures on the real
axis is still an open problem.

\end{document}